# $p$-TOPOLOGICAL AND $p$-REGULAR: DUAL NOTIONS IN CONVERGENCE THEORY

## SCOTT A. WILDE and D. C. KENT

ABSTRACT. The natural duality between "topological" and "regular," both considered as convergence space properties, extends naturally to $p$-regular convergence spaces, resulting in the new concept of a $p$-topological convergence space. Taking advantage of this duality, the behavior of $p$-topological and $p$-regular convergence spaces is explored, with particular emphasis on the former, since they have not been previously studied. Their study leads to the new notion of a neighborhood operator for filters, which in turn leads to an especially simple characterization of a topology in terms of convergence criteria. Applications include the topological and regularity series of a convergence space.

**Introduction.** In 1990, G. Richardson and one of the authors introduced the notion of *$p$-regular convergence space*, [6], defined as follows: If $q$ and $p$ are convergence structures on a set $X$, then the space $(X,q)$ is $p$-regular if $\text{cl}_p \mathcal{F} \xrightarrow{q} x$ whenever $\mathcal{F} \xrightarrow{q} x$, where "$\text{cl}_p$" is the $p$-closure operator. Clearly $p$-regularity is equivalent to regularity when $p = q$. By varying $p$, one can characterize various convergence properties in terms of $p$-regularity (see [6, 7]).

More recently, Kent and Richardson [7] developed some ideas and results due to Kowalsky [8], Cook and Fischer [1], and Biesterfeldt [2] to give convergence characterizations of the properties "topological" and "regular" so as to reveal a fundamental duality between these notions. These characterizations made use of "diagonal" axioms **F** and **R** which are in a natural way dual to each other. (It should be noted that the axiom called **R** in this paper was called **DF** in [7].)

In this paper, we begin by proving the $p$-regularity of a convergence space $(X,q)$ also has a "diagonal" characterization in terms of an axiom we call $\mathbf{R}_{p,q}$, which is obtained by making a minor alteration in the axiom **R**. We then use the dual axiom $\mathbf{F}_{p,q}$ to define (and introduce) the dual notion of a "$p$-topological convergence space."

Our goal is two-fold. We wish to study and develop this new concept of a $p$-topological convergence space, while simultaneously exploring the duality alluded to in the title of the paper. The approach based on duality is most useful in examining the structural behavior of $p$-topological and $p$-regular spaces as well as their upper and lower modifications. This approach is adopted in Sections 1 and 4. In Section 2, we study some aspects of $p$-topological spaces which do not have obvious analogues in the setting of $p$-regular spaces. Section 3 introduces the "neighborhood operator for filters" which seems to be "tailor-made" for the study of $p$-topological spaces and is used extensively in Section 4. The characterization of $p$-topological spaces, given in Theorem 3.2, yields a corollary which gives a simple and elegant characterization of a topology in terms of convergence criteria.

As is shown in Section 2 of this paper and also in [6, 7], both of the properties "$p$-topological" and "$p$-regular" can be adapted to characterize various convergence and topological concepts and, thereby, reveal underlying relationships between them. Other applications of these notions include the *regularity* and *topological series* of a convergence space which are discussed briefly in Section 5.

**1. The Axioms $F_{p,q}$ and $R_{p,q}$.** For standard notation and terminology pertaining to convergence spaces, the reader is referred to [7]. In particular, $\mathbf{F}(X)$ denotes the set of all filters on a set $X$, $\mathbf{U}(X)$ the set of all ultrafilters on $X$, and $\mathbb{C}(X)$ the complete lattice of all convergence structures on $X$ (with the discrete topology as the greatest element). Let $\dot{x}$ denote the fixed ultrafilter on $X$ generated by $x \in X$.

If $(X,q)$ is a convergence space and $J$ an arbitrary set, let $\mathcal{F} \in \mathbf{F}(J)$ and let $\sigma : J \to \mathbf{F}(X)$ be an arbitrary "selection function." We define $\kappa\sigma\mathcal{F}$ to be the filter $\cup_{F \in \mathcal{F}} \cap_{x \in F} \sigma(x)$ in $\mathbf{F}(X)$; $\kappa\sigma\mathcal{F}$ is called the *compression of $\mathcal{F}$ relative to $\sigma$*.

We, next, define two axioms pertaining to two convergence structures $p,q$ on a set $X$.

$\mathbf{F}_{p,q}$: Let $J$ be any set, $\psi: J \to X$, and let $\sigma : J \to \mathbf{F}(X)$ have the property that $\sigma(y) \xrightarrow{p} \psi(y)$ for all $y \in J$. If $\mathcal{F} \in \mathbf{F}(J)$ is such that $\psi(\mathcal{F}) \xrightarrow{q} x$, then $\kappa\sigma\mathcal{F} \xrightarrow{q} x$.

$\mathbf{R}_{p,q}$: Let $J$ be any set, $\psi: J \to X$ and let $\sigma : J \to \mathbf{F}(X)$ have the property that $\sigma(y) \xrightarrow{p} \psi(y)$, for all $y \in J$. If $\mathcal{F} \in \mathbf{F}(J)$ is such that $\kappa\sigma\mathcal{F} \xrightarrow{q} x$, then $\psi(\mathcal{F}) \xrightarrow{q} x$.

**THEOREM 1.1.** *Let $(X,q)$ be a convergence space and $p \in \mathbb{C}(X)$. Then $(X,q)$ is $p$-regular if and only if $p$ and $q$ satisfy $\mathbf{R}_{p,q}$.*

**PROOF.** Recall that $(X,q)$ is $p$-regular if $\text{cl}_p(\mathcal{F}) \xrightarrow{q} x$ whenever $\mathcal{F} \xrightarrow{q} x$.

($\Leftarrow$) Assume that $(X,q)$ and $p$ satisfy $\mathbf{R}_{p,q}$. Let $J = \{(\mathcal{G},y): \mathcal{G} \in \mathbf{U}(X), y \in X, \mathcal{G} \xrightarrow{p} y\}$. Define $\psi: J \to X$ by $\phi(\mathcal{G},y) = y$ and $\sigma: J \to \mathbf{F}(X)$ by $\sigma(\mathcal{G},y) = \mathcal{G}$. Note that $\sigma(z) \xrightarrow{p} \psi(z)$, for all $z = (\mathcal{G},y) \in J$.

Assume that $\mathcal{F} \xrightarrow{q} x$. We define a filter $\mathcal{H} \in \mathbf{F}(J)$ as follows: for each $F \in \mathcal{F}$, let $H_F = \{(\mathcal{G},y) \in J: F \in \mathcal{G}\}$, and let $\mathcal{H}$ be the filter on $J$ generated by $\{H_F: F \in \mathcal{F}\}$. Since $F \in \sigma(\mathcal{G},y)$ for every $(\mathcal{G},y) \in H_F$, $F \in \kappa\sigma\mathcal{H}$, so $\kappa\sigma\mathcal{H} \geq \mathcal{F}$. Thus, $\kappa\sigma\mathcal{H} \xrightarrow{q} x$. But observe that $\psi(H_F) = \text{cl}_p(F)$, so $\text{cl}_p(\mathcal{F}) \geq \psi(\mathcal{H})$. By $\mathbf{R}_{p,q}$, $\psi(\mathcal{H}) \xrightarrow{q} x$, which then implies that $\text{cl}_p(\mathcal{F}) \xrightarrow{q} x$. Thus, $(X,q)$ is $p$-regular.

($\Rightarrow$) Assume that $(X,q)$ is $p$-regular. Let $J, \sigma, \phi$ be as in $\mathbf{R}_{p,q}$ and let $\mathcal{F} \in \mathbf{F}(J)$ such that $\kappa\sigma\mathcal{F} \xrightarrow{q} x$. We claim that $\text{cl}_p(\kappa\sigma\mathcal{F}) \leq \psi(\mathcal{F})$. Let $F \in \mathcal{F}$ and choose $A_y \in \sigma(y)$, for every $y \in F$. Then $\sigma(y) \xrightarrow{p} \psi(y)$, for every $y \in F$, which implies that $\psi(y) \in \text{cl}_p(\cup_{y \in F} A_y)$ holds for every $y \in F$, and so $\psi(F) \subseteq \text{cl}_p(\cup_{y \in F} A_y)$. Since $\cup_{y \in F} A_y$ is a basic set in $\kappa\sigma\mathcal{F}$, the claim is verified. By $p$-regularity, $\text{cl}_p(\kappa\sigma\mathcal{F}) \xrightarrow{q} x$, which implies that $\psi(\mathcal{F}) \xrightarrow{q} x$. □

If $(X,q)$ is a convergence space and $p \in \mathbb{C}(X)$, then $(X,q)$ is defined to be $p$-

*topological* if $(X,q)$ and $p$ satisfy the axiom $\mathbf{F}_{p,q}$. Note that, by Theorem 1.1, $(X,q)$ is *p*-regular if and only if $(X,q)$ and $p$ satisfy $\mathbf{R}_{p,q}$. Since $\mathbf{F}_{p,q}$ and $\mathbf{R}_{p,q}$ are dual to each other, "*p*-topological" and "*p*-regular" are likewise dual properties. In the special case where $p = q$, $\mathbf{F}_{q,q}$ and $\mathbf{R}_{q,q}$ are denoted by $\mathbf{F}$ and $\mathbf{R}$, respectively.

**THEOREM 1.2.** *Let $(X,q)$ be a convergence space.*
  (i) *$(X,q)$ is topological if and only if $(X,q)$ satisfies $\mathbf{F}$.*
  (ii) *$(X,q)$ is regular if and only if $(X,q)$ satisfies $\mathbf{R}$.*

**PROOF.** The first assertion is proved in [6], the second by combining results from [2, 1]. □

It follows from Theorem 1.2 that "*p*-topological" generalizes "topological" in the same way that "*p*-regular" generalizes "regular". In the next theorem, $\mathbf{F}_{p,q}$ and $\mathbf{R}_{p,q}$ are applied directly to determine the behavior of these properties relative to initial constructs.

**THEOREM 1.3.** *Initial structures.*

   (i) *Let $\{(X_i,q_i): i \in I\}$ be a set of spaces together with a set of convergence structures $p_i$ which satisfies $\mathbf{F}_{p_i,q_i}$, for all $i \in I$. Let $X$ be a set and let $f_i: X \to X_i$ be a mapping, for each $i \in I$. If $q$ is the initial structure on $X$ relative to the families $\{(X_i,q_i): i \in I\}$ and $\{f_i: i \in I\}$, and $p$ is the initial structure on $X$ relative to $\{p_i: i \in I\}$ and $\{f_i: i \in I\}$, then $(X,q)$ and $p$ satisfy $\mathbf{F}_{p,q}$.*

   (ii) *Statement* (i) *remains valid if $\mathbf{F}_{p_i,q_i}$ is replaced by $\mathbf{R}_{p_i,q_i}$ and $\mathbf{F}_{p,q}$ is replaced by $\mathbf{R}_{p,q}$.*

**PROOF.** (i) It is well known that $q$-convergence is characterized by: $\mathcal{F} \xrightarrow{q} x$ if and only if $f_i(\mathcal{F}) \xrightarrow{q_i} f_i(x)$, for all $i \in I$. Let $J$ be a set and $\psi: J \to X$ and $\sigma: J \to \mathbf{F}(X)$ have the property that $\sigma(j) \xrightarrow{p} \psi(j)$ for all $j \in J$. Define $\sigma_i(j)$ and $\psi_i(j)$ so that $\sigma_i(j) = f_i(\sigma(j))$ and $\psi_i(j) = f_i(\psi(j))$ for all $j \in J$ and $i \in I$, respectively. Thus, $\sigma_i(j) \xrightarrow{p_i} \psi_i(j)$ for all $j \in J$. Also, $f_i(\kappa\sigma\mathcal{F}) = \kappa(f_i \circ \sigma)\mathcal{F} = \kappa\sigma_i\mathcal{F}$. Now, let $\mathcal{F} \in \mathbf{F}(J)$ have the property that $\psi(\mathcal{F}) \xrightarrow{q} x$ which then implies that $f_i(\psi(\mathcal{F})) \xrightarrow{q_i} f_i(x)$ for all $i \in I$, by the property of $q$ being the initial structure of all the $q_i$. Thus, $f_i(\kappa\sigma\mathcal{F}) = \kappa\sigma_i\mathcal{F} \xrightarrow{q_i} f_i(x)$ for all $i \in I$ by the property $\mathbf{F}_{p_i,q_i}$. Hence, $\kappa\sigma\mathcal{F} \xrightarrow{q} x$ by the definition of $q$ and this implies that $(X,q)$ and $p$ satisfy $\mathbf{F}_{p,q}$.

  (ii) This proof is essentially the same as that of (i). □

**COROLLARY 1.4.** *A subspace of a p-topological (respectively, p-regular) space is $p'$-topological (respectively, $p'$-regsular), where $p'$ denotes the restriction of $p$ to the subspace.*

**COROLLARY 1.5.** *Let $(X,q) = \Pi_{i \in I}(X_i,q_i)$ and $(X,p) = \Pi_{i \in I}(X_i,p_i)$ be product convergence spaces. If each $(X_i,q_i)$ is $p_i$-topological (respectively, $p_i$-regular), then $(X,q)$ is p-topological (respectively, p-regular).*

**COROLLARY 1.6.** *Let $X$ be a set and let $\Lambda = \{q_i: i \in I\}$ and $\Gamma = \{p_i: i \in I\}$ be subsets of $\mathbb{C}(X)$. Let $q = \sup\Lambda$ and $p = \sup\Gamma$. If $(X_i,q_i)$ is $p_i$-topological (respectively, $p_i$-regular) for each $i \in I$, then $(X,q)$ is p-topological (respectively, p-regular).*

Before proving the analogue of Theorem 1.3 for final structures, we give a simpler characterization for $p$-topological spaces which makes use of the $p$-interior operator $\mathcal{I}_p$.

**THEOREM 1.7.** *Let $(X,q)$ be a convergence space and $p \in \mathbb{C}(X)$. Then $(X,q)$ is $p$-topological if and only if, whenever $\mathcal{F} \xrightarrow{q} x$, there exists $\mathcal{G} \xrightarrow{q} x$ such that $\mathcal{F} \geq \mathcal{I}_p\mathcal{G}$.*

**PROOF.** Assume that $(X,q)$ and $p$ satisfy $\mathbf{F}_{p,q}$. Let $\mathcal{F} \xrightarrow{q} x$ and $J = \{(\mathcal{G},x), : \mathcal{G} \in \mathbf{U}(X), \mathcal{G} \xrightarrow{p} x\}$. Let $\psi : J \to X$ be defined by $\psi((\mathcal{G},x)) = x$ and $\sigma : J \to \mathbf{F}(X)$ by $\sigma((\mathcal{G},x)) = \mathcal{G}$. Note that $\psi$ is onto $X$ (since $x \in X \Rightarrow (\dot{x},x) \in J$). If $\mathcal{H} = \psi^{-1}(\mathcal{F})$, then $\mathcal{H} \in \mathbf{F}(J)$ and $\psi(\mathcal{H}) = \mathcal{F} \xrightarrow{q} x$. By $\mathbf{F}_{p,q}$, $\kappa\sigma\mathcal{H} \xrightarrow{q} x$. Thus, to show that $(X,q)$ satisfies the given condition, it suffices to show that $\mathcal{F} \geq \mathcal{I}_p(\kappa\sigma\mathcal{H})$.

Let $F \in \mathcal{F}$; then $\psi^{-1}(F)$ is a basic set in $\mathcal{H}$. Note that $\mathcal{G} \in \mathbf{U}(X)$ and $\mathcal{G} \xrightarrow{p} x \in F$ imply that $(\mathcal{G},x) \in \psi^{-1}(F)$. Since $\sigma((\mathcal{G},x)) = \mathcal{G}$, for each pair $(\mathcal{G},x) \in \psi^{-1}(F)$, choose $G_{(\mathcal{G},x)} \in \mathcal{G}$. Then $A = \cup\{G_{(\mathcal{G},x)} : (\mathcal{G},x) \in \psi^{-1}(F)\}$ is a basic set in $\kappa\sigma\mathcal{H}$. For a given $y \in F$, $A_y = \cup\{G_{(\mathcal{G},y)} : \mathcal{G} \in \mathbf{U}(X), \mathcal{G} \xrightarrow{p} y\} \in \mathcal{V}_p(y)$ (since $\mathcal{V}_p(y)$ is the intersection of all ultrafilters which $p$-converge to y). Since $A_x \subseteq A$, for all $x \in F$, $A \in \mathcal{V}_p(x)$, for all $x \in F$. Thus, $F \subseteq \mathcal{I}_p(A)$ and we obtain the desired conclusion that $\mathcal{F} \geq \mathcal{I}_p(\kappa\sigma\mathcal{H})$.

Conversely, let $J, \psi, \sigma$, and $\mathcal{F}$ be as in $\mathbf{F}_{p,q}$ and let $\psi(\mathcal{F}) \xrightarrow{q} x$. Since $(X,q)$ satisfies the specified condition, there exists a filter $\mathcal{G} \xrightarrow{q} x$ such that $\psi(\mathcal{F}) \geq \mathcal{I}_p\mathcal{G}$. To show that $\kappa\sigma\mathcal{F} \xrightarrow{q} x$, it suffices to show that $\kappa\sigma\mathcal{F} \geq \mathcal{G}$. Let $G \in \mathcal{G}$ and choose $F \in \mathcal{F}$ such that $\psi(F) \subseteq \mathcal{I}_p(G)$. For each $y \in F$, $\psi(y) \in \mathcal{I}_p(G)$ implies that $G \in \sigma(y)$. Thus, $G \in \kappa\sigma\mathcal{F}$, which yields the desired result that $\kappa\sigma\mathcal{F} \geq \mathcal{G}$. □

Let $f : (X,q) \to (Y,p)$ be a function between convergence spaces. We define $f$ to be an *interior map* if $f(\mathcal{I}_p(A)) \subseteq \mathcal{I}_p(f(A))$ holds for all $A \subseteq X$, and a *closure map* if $\mathrm{cl}_p(f(A)) \subseteq f(\mathrm{cl}_q(A))$ holds for all $A \subseteq X$. Closure maps were introduced in [6], where they were found to be useful in the study of $p$-regularity.

**THEOREM 1.8.** *Let $X$ be a set, $\{(X_i,q_i) : i \in I\}$ a set of convergence spaces, and $\{f_i : i \in I\}$ a set of functions mapping $X_i$ to $X$ such that $X = \cup_{i \in I} f(X_i)$. Let $q$ be the final convergence structure on $X$ induced by $\{f_i : i \in I\}$ and $\{(X_i,q_i) : i \in I\}$.*

*(i) If each $(X_i,q_i)$ is $p_i$-topological for some $p_i \in \mathbb{C}(X_i)$ and $p$ is a convergence structure on $X$ such that each $f_i : (X_i,p_i) \to (X,p)$ is an interior map, then $(X,q)$ is $p$-topological.*

*(ii) If each $(X_i,q_i)$ is $p_i$-regular for some $p_i \in \mathbb{C}(X_i)$ and $p$ is a convergence structure on $X$ such that each $f_i : (X_i,p_i) \to (X,p)$ is a closure map, then $(X,q)$ is $p$-regular.*

**PROOF.** (i) Let $\mathcal{F} \xrightarrow{q} x$. Then there exists $j \in I$, $x_j \in X_j$ such that $f_j(x_j) = x$ and $\mathcal{F}_j \xrightarrow{q_j} x_j$ such that $f_j(\mathcal{F}_j) \leq \mathcal{F}$. Since $(X_j,q_j)$ is $p_j$-topological, there exists $\mathcal{G}_j \xrightarrow{q_j} x_j$ such that $\mathcal{F}_j \geq \mathcal{I}_{p_j}(\mathcal{G}_j)$. Since $f_j$ is an interior map, $f_j(\mathcal{I}_{p_j}(\mathcal{G}_j)) \geq \mathcal{I}_p(f_j(\mathcal{G}_j))$. By continuity of $f_j : (X_j,q_j) \to (X,q)$, $f_j(\mathcal{G}_j) \xrightarrow{q} x$, and $\mathcal{F} \geq f_j(\mathcal{F}_j) \geq \mathcal{I}_p(f_j(\mathcal{G}_j))$, so $(X,q)$ is $p$-topological by Theorem 1.7.

The proof of (ii) is similar. □

**COROLLARY 1.9.** *Let $f : (X',q') \to (X,q)$ be a convergence quotient map.*
*(i) If $f : (X',p') \to (X,p)$ is an interior map and $(X',q')$ is $p$-topological, then $(X,q)$*

is *p*-topological.

(ii) *If $f: (X', p') \longrightarrow (X, p)$ is a closure map and $(X', q')$ is p-regular, then $(X, q)$ is p-regular.*

**COROLLARY 1.10.** *Let $(X, q) = \sum_{i \in I}(X_i, q_i)$ be a disjoint sum of convergence spaces and $p \in \mathbb{C}(X)$. Let $g_i: X_i \longrightarrow X$ be the canonical injection.*

(i) *If for each $i \in I$, $(X_i, q_i)$ is $p_i$-topological and $g_i: (X_i, p_i) \longrightarrow (X, p)$ is an interior map, then $(X, q)$ is p-topological.*

(ii) *If for each $i \in I$, $(X_i, q_i)$ is $p_i$-regular and $g_i: (X_i, p_i) \longrightarrow (X, p)$ is a closure map, then $(X, q)$ is p-regular.*

**COROLLARY 1.11.** *Let $\Lambda = \{q_i: i \in I\} \subseteq \mathbb{C}(X)$, let $p \in \mathbb{C}(X)$, and assume that $(X, q_i)$ is $p_i$-topological (respectively, $p_i$-regular), for all $i \in I$. If $p_i \leq p$ for each $i \in I$ and $q = \inf \Lambda$, then $(X, q)$ is p-topological (respectively, p-regular).*

The final result of this section, which follows immediately from Corollaries 1.6 and 1.11, asserts that for a fixed convergence structure $p$ on $X$, both of the properties "*p*-topological" and "*p*-regular" are preserved under arbitrary infima and suprema in the lattice $\mathbb{C}(X)$.

**COROLLARY 1.12.** *Let $\Lambda = \{q_i: i \in I\} \subseteq \mathbb{C}(X)$ and let $p \in \mathbb{C}(X)$ be such that $(X_i, q_i)$ is p-topological (respectively, p-regular), for all $i \in I$. Let $q = \inf \Lambda$ and $r = \sup \Lambda$. Then both $(X, q)$ and $(X, r)$ are p-topological (respectively, p-regular).*

**2. More on *p*-topological spaces.** In Section 1, we observed that *p*-topological and *p*-regular properties exhibit essentially the same structural behavior. Now, we gain some additional insight into the behavior of *p*-topological spaces by making use of Theorem 1.7. The first result of this section gives a simple characterization of pretopological spaces which are *p*-topological.

**THEOREM 2.1.** *Let $(X, q)$ be a pretopological space and $p \in \mathbb{C}(X)$.*

(i) *$(X, q)$ is p-topological if and only if $\mathcal{V}_q(x) = \mathcal{I}_p \mathcal{V}_q(x)$.*

(ii) *If $(X, q)$ is p-topological, then $q \leq \tau p$, where $\tau p$ denotes the topological modification of $p$.*

**PROOF.** (i) Assume that $(X, q)$ is *p*-topological. Since $\mathcal{V}_q(x) \xrightarrow{q} x$, it follows by Theorem 1.7 that $\mathcal{V}_q(x) \geq \mathcal{I}_p \mathcal{V}_q(x)$, and, hence, $\mathcal{V}_q(x) = \mathcal{I}_p \mathcal{V}_q(x)$. Conversely, if the given equality holds, then $\mathcal{F} \xrightarrow{q} x$ implies $\mathcal{F} \geq \mathcal{V}_q(x) = \mathcal{I}_p \mathcal{V}_q(x)$, and since $(X, q)$ is pretopological, $\mathcal{V}_q(x) \xrightarrow{q} x$, so $(X, q)$ is *p*-topological by Theorem 1.7.

(ii) If $(X, q)$ is *p*-topological, then by (i) $\mathcal{V}_q(x) = \mathcal{I}_p \mathcal{V}_q(x)$, and it follows that $\mathcal{V}_q(x)$ has a filter base of *p*-open sets (which are the same as $\tau p$-open sets). Thus, $\mathcal{V}_q(x) \leq \mathcal{V}_{\tau p}(x)$, and since $q$ is a pretopology, $q \leq \tau p$. □

**EXAMPLE 2.2.** Converse of Theorem 2.1(ii) is generally false.

Let $X = \mathbb{R}$ be the set of real numbers, and let $\tau$ denote the usual topology on $\mathbb{R}$. Note that $\tau \cup \{0\}$ is a base for a topology $p$ on $\mathbb{R}$, where $\tau < p$ and $\mathcal{V}_p(x) = \mathcal{V}_\tau(x)$, for all $x \neq 0$, whereas $\mathcal{V}_p(0) = \dot{0}$. Let $q$ be the pretopology on $\mathbb{R}$ defined by $\mathcal{V}_q(x) = \mathcal{V}_\tau(x)$ for $x \neq 0$ and $\mathcal{V}_q(0) = \mathcal{V}_\tau(0) \vee \dot{Q}$ (where $\dot{Q}$ is the filter of oversets of the set $Q$ of rational numbers). Note that $\tau < q < p$. Then $\mathcal{V}_q(x) \xrightarrow{q} x$, but $\mathcal{I}_p \mathcal{V}_q(0) = \dot{0} \neq \mathcal{V}_q(0)$, so by Theorem 2.1(i), $(X,q)$ is not $p$-topological.

**COROLLARY 2.3.** *If $p$ and $q$ are topological, then $(X,q)$ is $p$-topological if and only if $q \leq p$.*

**PROOF.** If $(X,q)$ is $p$-topological, then $p \leq q$ follows from Theorem 2.1(ii). Conversely, if $q \leq p$, then $\mathcal{I}_p \mathcal{V}_q(x) = \mathcal{V}_q(x)$ follows because $q$ is a topology, and so, the conclusion follows from Theorem 2.1(i). □

The preceding example shows that Corollary 2.3 does not hold under the weaker condition that $q$ is pretopological.

Note that if $(X,q)$ is $p$-topological, then $(X,q)$ is obviously $p'$-topological for any $p' \geq p$. Clearly, every convergence space is $\delta$-topological, where $\delta$ denotes the discrete topology.

**COROLLARY 2.4.** *If $(X,q)$ is $p$-topological, then $(X,\pi q)$ and $(X,\tau q)$ are $p$-topological, and $\tau q \leq \pi q \leq \tau p$ (where $\pi q$ denotes the pretopological modification of $q$).*

**PROOF.** Let $\mathcal{F} \xrightarrow{q} x$; then, by Theorem 1.7, there exists $\mathcal{G} \xrightarrow{q} x$ such that $\mathcal{F} \geq \mathcal{I}_p \mathcal{G} \geq \mathcal{I}_p \mathcal{V}_q(x)$. This holds for every $\mathcal{F} \xrightarrow{q} x$, so $\mathcal{V}_{\pi q}(x) = \mathcal{V}_q(x) \geq \mathcal{I}_p \mathcal{V}_{\pi q}(x)$. Thus, $\pi q$ is $p$-topological from Theorem 2.1(i). $\tau q \leq \pi q \leq \tau p$ follows from Theorem 2.1(ii), and $(X,\tau q)$ is $\tau p$-topological from Corollary 2.3, $(X,\tau q)$ is $p$-topological from the remark preceding the corollary, since $\tau p \leq p$. □

**COROLLARY 2.5.** *Let $\iota$ denote the indiscrete topology on $X$. Then $(X,\iota)$ is $p$-topological, for every $p \in \mathbb{C}(X)$.*

**PROOF.** By Corollary 2.3, $(X,\iota)$ is $\iota$-topological, and, hence, $p$-topological for all $p \in \mathbb{C}(X)$ by the remark preceding Corollary 2.4. □

Given a convergence space $(X,q)$, let $\rho q$ denote the finest completely regular topology on $X$ coarser than $q$, and let $\omega q$ be the finest completely regular topology on $X$ coarser than $q$.

**THEOREM 2.6.** *A convergence space $(X,q)$ is a regular (respectively, completely regular) topological space if and only if $(X,q)$ is $\rho q$-topological (respectively, $\omega q$-topological).*

**PROOF.** If $(X,q)$ is a regular topological space, then $q = \rho q$ and $(X,q)$ is obviously $q$-topological. Hence, $\rho q$-topological. Conversely, if $(X,q)$ is $\rho q$-topological, then $(X,q)$ is clearly $q$-topological, and, hence, topological. By Corollary 2.4, $q \leq \rho q$, and, hence, $q = \rho q$. Thus, $(X,q)$ is regular and topological. □

Let $(X,q)$ be a topological space, and let $q'$ be the topology on $X$ generated by $\mathcal{B}_{q'} = \{X\} \cup \{U \subseteq X : U \in q \text{ and } U \subseteq K \text{ for some } q\text{-compact subset } K \text{ of } X\}$.

**THEOREM 2.7.** *A $T_1$ topological space $(X,q)$ is locally compact if and only if $(X,q)$ is $q'$-topological.*

**PROOF.** Let $(X,q)$ be locally compact and $x \in X$. Let $U \in \mathcal{V}_q(x)$ be $q$-open. By local compactness, there is a compact set $A \in \mathcal{V}_q(x)$. Let $V$ be a $q$-open set such that $V \subseteq U \cap A$. Then $V$ is $q'$-open. So, $\mathcal{I}_{q'}\mathcal{V}_q(x) = \mathcal{V}_q(x)$, which implies, by Theorem 2.1(i), that $(X,q)$ is $q'$-topological.

Conversely, let $(X,q)$ be $q'$-topological and $x \in X$. Since $q$ is $T_1$, there exists $U \in \mathcal{V}_q(x)$ such that $U \ne X$. Since $\mathcal{I}_{q'}\mathcal{V}_q(x) = \mathcal{V}_q(x)$, by Theorem 2.1(i), $\mathcal{I}_{q'}(U) \in \mathcal{V}_q(x)$, and $\mathcal{I}_{q'}(U) = U \subseteq A$, which implies that $A \in \mathcal{V}_q(x)$. Thus, $(X,q)$ is locally compact. □

A related theorem characterizing local compactness in terms of $p$-regularity is the following result, which is a direct corollary of [6, Thm. 3.1].

**THEOREM 2.8.** *Let $(X,q)$ be a $T_1$ convergence space. Let $p$ be the topology on $X$ having as a base of closed sets all the nonempty subsets of $q$-compact sets. Then $(X,q)$ is locally compact if and only if $(X,q)$ is $p$-regular.*

**3. The neighborhood operator for a filter.** In this section, we introduce a new filter notion which is essentially dual to the "closure of a filter," thereby obtaining another characterization of "$p$-topological" which further illustrates its duality with "$p$-regular".

Let $(X,q)$ be a convergence space, $\mathcal{F} \in \mathbf{F}(X)$. Then $\mathcal{V}_q\mathcal{F} = \{A \in \mathcal{F}: \mathcal{I}_q(A) \in \mathcal{F}\}$ is called the *$q$-neighborhood filter* of $\mathcal{F}$.

**PROPOSITION 3.1.** *If $(X,q)$ is a convergence space and $\mathcal{F} \in \mathbf{F}(X)$, then $\mathcal{V}_q\mathcal{F}$ is the finest filter on $X$ such that $\mathcal{F} \ge \mathcal{I}_q(\mathcal{V}_q\mathcal{F})$.*

**PROOF.** It is clear that $\mathcal{V}_q\mathcal{F}$ is a filter on $X$ such that $\mathcal{I}_q\mathcal{V}_q\mathcal{F} \le \mathcal{F}$. If $\mathcal{G}$ is any filter on $X$ such that $\mathcal{I}_q\mathcal{G} \le \mathcal{F}$, then $G \in \mathcal{G}$ implies $\mathcal{I}_q G \in \mathcal{F}$, and, hence, $G \in \mathcal{V}_q\mathcal{F}$. □

If $\mathcal{F} = \dot{x}$, it is obvious from the definition that $\mathcal{V}_q\dot{x} = \mathcal{V}_q(x)$ is the $q$-neighborhood filter at $x$.

Recall that $(X,q)$ is $p$-regular if $\mathcal{F} \xrightarrow{q} x$ implies $\mathrm{cl}_p \mathcal{F} \xrightarrow{q} x$. The corresponding dual characterization for a $p$-topological space is the following.

**THEOREM 3.2.** *A convergence space $(X,q)$ is $p$-topological if and only if $\mathcal{F} \xrightarrow{q} x$ implies $\mathcal{V}_p\mathcal{F} \xrightarrow{q} x$.*

**PROOF.** Let $(X,q)$ be $p$-topological and $\mathcal{F} \xrightarrow{q} x$. By Theorem 1.7, there is $\mathcal{G} \xrightarrow{q} x$ such that $\mathcal{F} \ge \mathcal{I}_q\mathcal{G}$. By Proposition 2.1, $\mathcal{V}_p\mathcal{F} \ge \mathcal{G}$, and so $\mathcal{V}_p\mathcal{F} \xrightarrow{q} x$. Conversely if the condition holds, we can set $\mathcal{G} = \mathcal{V}_p\mathcal{F}$ in Theorem 1.7, and, thus, $(X,q)$ is $p$-topological. □

**COROLLARY 3.3.** *A convergence space $(X,q)$ is topological if and only if $\mathcal{F} \xrightarrow{q} x$ implies $\mathcal{V}_q\mathcal{F} \xrightarrow{q} x$.*

The $q$-neighborhood filter of a filter can also be described by means of the compression operator for filters defined in Section 1.

**PROPOSITION 3.4.** *Let $(X,p)$ be a convergence space and let $\sigma: X \to \mathbf{F}(X)$ be defined by $\sigma(x) = \mathcal{V}_p(x)$ for all $x \in X$. Then for any $\mathcal{F} \in \mathbf{F}(X)$, $\kappa\sigma\mathcal{F} = \mathcal{V}_p(\mathcal{F})$.*

**PROOF.** Let $A \in \mathcal{V}_p(\mathcal{F})$. Then $\mathcal{I}_p(A) \in \mathcal{F}$. If $F = \mathcal{I}_p(A)$, then for each $x \in F$, $A \in$ $\mathcal{V}_p(x)$, and so $A = \bigcup_{x \in F} V_x$, where each $V_x = A$, is a basic set in $\kappa\sigma\mathcal{F}$. Conversely, let $A \in \kappa\sigma\mathcal{F}$. Then $A$ contains a basic set of the form $B = \bigcup_{y \in F} V_y$, where $F \in \mathcal{F}$ and $V_y \in \mathcal{V}_p(y)$, for all $y \in F$. To show that $A \in \mathcal{V}_p(\mathcal{F})$, it suffices to show that $F \subseteq \mathcal{I}_p(B)$. $x \in F$ implies that $x \in V_x \subseteq B$. Thus, $B \in \mathcal{V}_p(x)$ and so, $x \in \mathcal{I}_p(B)$. □

Let $(X,q)$ be a convergence space and $\mathcal{F} \in \mathbf{F}(X)$. For any $n \in \mathbb{N}$, the set of natural numbers, the $n$th iterations of the closure and neighborhood operators for a filter $\mathcal{F}$ are given inductively by:

$$\begin{aligned} \mathrm{cl}_q^n \mathcal{F} &= \mathrm{cl}_q\left(\mathrm{cl}_q^{n-1}\mathcal{F}\right), \\ \mathcal{V}_q^n \mathcal{F} &= \mathcal{V}_q\left(\mathcal{V}_q^{n-1}\mathcal{F}\right). \end{aligned} \qquad (1)$$

The next two propositions summarize (without proof) some additional elementary properties of the neighborhood operator for filters.

**PROPOSITION 3.5.** *Let $(X,q)$ be a convergence space, $n \in \mathbb{N}$, and $\{\mathcal{F}_i : i \in I\} \subseteq \mathbf{F}(X)$. Then:*

(i) $\mathcal{V}_p^n(\cap_{i \in I} \mathcal{F}_i) = \cap_{i \in I} \mathcal{V}_p^n \mathcal{F}_i$;

(ii) *If $\vee_{i \in I} \mathcal{F}_i$ exists, then $\mathcal{V}_p^n(\vee_{i \in I} \mathcal{F}_i) \geq \vee_{i \in I} \mathcal{V}_p^n \mathcal{F}_i$;*

(iii) *Equality holds in (ii) under the additional assumption that $\{\mathcal{F}_i : i \in I\}$ is an upward directed set of filters.*

**PROPOSITION 3.6.** *Let $f : (X,q) \to (Y,p)$ be a function between convergence spaces. Let $\mathcal{F} \in \mathbf{F}(X)$ and $n \in \mathbb{N}$.*

(i) *If $f$ is continuous, then $f(\mathcal{V}_q^n \mathcal{F}) \geq \mathcal{V}_p^n f(\mathcal{F})$.*

(ii) *If $f$ is an interior map, then $f(\mathcal{V}_q^n \mathcal{F}) \leq \mathcal{V}_p^n f(\mathcal{F})$.*

**4. Lower and upper modifications.** It was established in Corollary 1.12 that each of the properties $p$-topological and $p$-regular is preserved under both infima and suprema in the lattice $\mathbb{C}(X)$. Since an indiscrete space is both $p$-topological and $p$-regular for any choice of $p$, we immediately obtain the following.

**PROPOSITION 4.1.** *Let $(X,q)$ be a convergence space and $p \in \mathbb{C}(X)$.*

(i) *There is a finest $p$-topological convergence structure $\tau_p q$ on $X$ coarser than $q$.*

(ii) *There is a finest $p$-regular convergence structure $r_p q$ on $X$ coarser than $q$.*

The structures $\tau_p q$ and $r_p q$ are called the *lower $p$-topological* and *lower $p$-regular modifications* of $q$, respectively. The dual relationship between these concepts is evident in the next theorem.

**THEOREM 4.2.** *Let $(X,q)$ be a convergence space and $p \in \mathbb{C}(X)$.*

(i) $\mathcal{F} \xrightarrow{\tau_p q} x$ *if and only if there exists* $\mathcal{G} \xrightarrow{q} x$ *such that* $\mathcal{F} \geq \mathcal{V}_p^n \mathcal{G}$, *for some* $n \in \mathbb{N}$.

(ii) $\mathcal{F} \xrightarrow{r_p q} x$ *if and only if there exists* $\mathcal{G} \xrightarrow{q} x$ *such that* $\mathcal{F} \geq \mathrm{cl}_p^n \mathcal{G}$, *for some* $n \in \mathbb{N}$.

**PROOF.** (i) Let $q'$ be defined by $\mathcal{F} \xrightarrow{q'} x$ if and only if there is $\mathcal{G} \xrightarrow{q} x$ such that $\mathcal{F} \geq \mathcal{V}_p^n(\mathcal{G})$, for some $n \in \mathbb{N}$. One may easily verify that $q'$ is a convergence structure. If $\mathcal{F} \xrightarrow{q} x$, then $\mathcal{F} \geq \mathcal{V}_p^n(\mathcal{F})$ for any $n \in \mathbb{N}$, and so $\mathcal{F} \xrightarrow{q'} x$. Thus, $q' \leq q$. To show that $q'$ is *p-topological*, let $\mathcal{F} \xrightarrow{q'} x$ and let $\mathcal{G} \xrightarrow{q} x$ be such that $\mathcal{F} \geq \mathcal{V}_p^n(\mathcal{G})$ for some $n \in \mathbb{N}$. Then $\mathcal{V}_p(\mathcal{F}) \geq \mathcal{V}_p^{n+1}(\mathcal{G})$, so $\mathcal{V}_p(\mathcal{F}) \xrightarrow{q} x$, and by Theorem 3.2, $q'$ is *p-topological*.

Finally, assume that $r$ is *p*-topological and $r \leq q$. Let $\mathcal{F} \xrightarrow{q'} x$. Then there is $\mathcal{G} \xrightarrow{q} x$ such that $\mathcal{F} \geq \mathcal{V}_p^n(\mathcal{G})$ for some $n \in N$. $\mathcal{G} \xrightarrow{q} x$ implies $\mathcal{G} \xrightarrow{r} x$, and since $r$ is *p-topological*, $\mathcal{V}_p^n(\mathcal{G}) \xrightarrow{r} x$, for all $n \in \mathbb{N}$. But $\mathcal{F} \geq \mathcal{V}_p^n(\mathcal{G})$ for some $n \in \mathbb{N}$, and, hence, $\mathcal{F} \xrightarrow{r} x$. Thus, $r \leq q'$ and the proof is complete.

(ii) See [6, Thm. 2.2]. □

Since the discrete topology $\delta$ on a set $X$ is generally neither *p*-topological nor *p*-regular for an arbitrary $p \in \mathbb{C}(X)$, the existence of an upper *p*-topological (or upper *p*-regular) modification for some $q \in \mathbb{C}(X)$ depends on the existence of a *p*-topological (or *p*-regular) convergence structure on $X$ finer than $q$. Clearly, $\tau_p \delta$ is the finest *p*-topological structure in $\mathbb{C}(X)$ and $r_p \delta$ is the finest *p*-regular member of $\mathbb{C}(X)$. Thus, a coarsest *p*-topological (respectively, *p*-regular) convergence structure on $X$ finer than $q$ exists if and only if $q \leq \tau_p \delta$ (respectively, $q \leq r_p \delta$). Using Theorem 4.2, this result may be restated as follows.

**THEOREM 4.3.** *Let $(X,q)$ be a convergence space and $p \in \mathbb{C}(X)$.*

(i) *There is a coarsest p-topological convergence structure $\tau^p q$ on $X$ finer than $q$ if and only if $\mathcal{V}_p^n(\dot{x}) \xrightarrow{q} x$, for all $x \in X$ and for all $n \in \mathbb{N}$.*

(ii) *There is a coarsest p-regular convergence structure $r^p q$ on $X$ finer than $q$ if and only if $\mathrm{cl}_p^n(\dot{x}) \xrightarrow{q} x$, for all $x \in X$ and for all $n \in \mathbb{N}$.*

When they exist, $\tau^p q$ and $r^p q$ are called the *upper p-topological* and *upper p-regular modifications* of $q$, respectively. Note that for $\tau^p q$ to exist, it is necessary that $q \leq p$, and that $r^p q$ will exist whenever $p$ is $T_1$.

**THEOREM 4.4.** *Let $(X,q)$ be a convergence space and $p \in \mathbb{C}(X)$.*

(i) *If $\tau^p q$ exists, then $\mathcal{F} \xrightarrow{\tau^p q} x$ if and only if $\mathcal{V}_p^n(\mathcal{F} \cap \dot{x}) \xrightarrow{q} x$, for all $n \in \mathbb{N}$.*

(ii) *If $r^p q$ exists, then $\mathcal{F} \xrightarrow{r^p q} x$ if and only if $\mathrm{cl}_p^n(\mathcal{F} \cap \dot{x}) \xrightarrow{q} x$, for all $n \in \mathbb{N}$.*

**PROOF.** (i) Let $q^*$ be defined by $\mathcal{F} \xrightarrow{q^*} x$ if and only if $\mathcal{V}_p^n(\mathcal{F} \cap \dot{x}) \xrightarrow{q} x$, for all $n \in \mathbb{N}$. It is easily shown that $q^*$ is a convergence structure. If $\mathcal{F} \xrightarrow{q^*} x$, then $\mathcal{F} \geq \mathcal{V}_p^n(\mathcal{F} \cap \dot{x})$, and $\mathcal{V}_p^n(\mathcal{F} \cap \dot{x}) \xrightarrow{q} x$ implies $\mathcal{F} \xrightarrow{q} x$. Thus, $q \leq q^*$. To show that $q^*$ is *p*-topological, $\mathcal{F} \xrightarrow{q^*} x$ implies $\mathcal{V}_p^{n+1}(\mathcal{F} \cap \dot{x}) = \mathcal{V}_p^n(\mathcal{V}_p(\mathcal{F} \cap \dot{x})) \xrightarrow{q} x$ for all $n \in \mathbb{N}$, which implies $\mathcal{V}_p^n(\mathcal{V}_p(x) \cap \dot{x}) \xrightarrow{q} x$, for all $n \in \mathbb{N}$. Thus, $\mathcal{V}_p(x) \xrightarrow{q^*} x$ and so, $q^*$ is *p-topological* by Theorem 3.2.

Finally, assume that $r$ is *p-topological* and $q \leq r$. Then, by Theorem 3.2, $\mathcal{F} \xrightarrow{r} x$ implies $\mathcal{V}_p^n(\mathcal{F} \cap \dot{x}) \xrightarrow{r} x$, for each $n \in \mathbb{N}$, and, hence, $\mathcal{V}_p^n(\mathcal{F} \cap \dot{x}) \xrightarrow{q} x$, for each $n \in \mathbb{N}$. But this implies $q^* \leq r$. Thus, $q^* = \tau^p q$.

The proof of (ii) is similar. □

**THEOREM 4.5.** *Let $(X,q)$ and $(X',q')$ be convergence spaces and let $f: (X,q) \longrightarrow (X',q')$ be continuous. Assume that $p \in \mathbb{C}(X)$ and $p' \in \mathbb{C}(X')$.*

(i) *If $f: (X,p) \longrightarrow (X',p')$ is continuous, then both of the mappings $f: (X, \tau_p q) \longrightarrow (X', \tau_{p'} q')$ and $f: (X, r_p q) \longrightarrow (X', r_{p'} q')$ are continuous.*

(ii) *If both of $\tau^p q$ and $\tau^{p'} q'$ exist and $f:(X,p) \longrightarrow (X',p')$ is an interior map, then $f:(X,\tau^p q) \longrightarrow (X,\tau^{p'} q')$ is continuous.*

(iii) *If both of $r^p q$ and $r^{p'} q'$ exist and $f:(X,p) \longrightarrow (X',p')$ is a closure map, then $f:(X,r^p q) \longrightarrow (X,r^{p'} q')$ is continuous.*

**PROOF.** Those results pertaining to $p$-regular structures have been proved in [6]. Those pertaining to $p$-topological structures can be proved analogously using Theorems 4.2(i) and 4.4(i), along with Proposition 3.6. □

The next two theorems show that the lower modifications behave reasonably well relative to final structures, whereas the upper modifications exhibit comparable behavior relative to initial structures.

**THEOREM 4.6.** *Let $X$ be a set and let $\{(X_i,q_i): i \in I\}$ and $\{(X_i,p_i): i \in I\}$ be collections of convergence spaces, and for all $i \in I$, $f_i: X_i \longrightarrow X$. Let $q$ be the final structure on $X$ induced by $\{f_i: i \in I\}$ and $\{(X_i,q_i): i \in I\}$ and let $p \in \mathbb{C}(X)$. Furthermore, assume that $X = \cup_{i\in I} f_i(X_i)$.*

(i) *If each $f_i: (X_i,p_i) \longrightarrow (X,p)$ is a continuous interior map, then $\tau_p q$ is the final structure on $X$ induced by $\{f_i: i \in I\}$ and $\{(X_i, \tau_{p_i} q_i): i \in I\}$.*

(ii) *If each $f_i: (X_i,p_i) \longrightarrow (X,p)$ is a continuous closure map, then $r_p q$ is the final structure on $X$ induced by $\{f_i: i \in I\}$ and $\{(X_i, r_{p_i} q_i): i \in I\}$.*

**PROOF.** (i) Let $s$ denote the final structure on $X$ induced by $\{f_i: i \in I\}$ and $\{(X_i, \tau_{p_i} q_i) : i \in I\}$. Let $\mathcal{F} \xrightarrow{s} x$. Then there is $i \in I$, $x_i \in X_i$, and $\mathcal{G}_i \xrightarrow{\tau_{p_i} q_i} x_i$ such that $\mathcal{F} \geq f_i(\mathcal{G}_i)$. Thus, there is $\mathcal{H}_i \xrightarrow{q_i} x_i$ and $n \in \mathbb{N}$ such that $\mathcal{G}_i \geq \mathcal{V}_{p_i}^n \mathcal{H}_i$ by Theorem 4.2. Hence, $f_i(\mathcal{H}_i) \xrightarrow{q} x$ and $\mathcal{F} \geq f_i(\mathcal{G}_i) \geq f_i(\mathcal{V}_{p_i}^n \mathcal{H}_i) = \mathcal{V}_p^n f_i(\mathcal{H}_i)$, where the last inequality follows by Proposition 3.6. Thus, $\mathcal{F} \xrightarrow{\tau_p q} x$.

Conversely, let $\mathcal{F} \xrightarrow{\tau_p q} x$. Then there is $\mathcal{G} \xrightarrow{q} x$ and $n \in \mathbb{N}$ such that $\mathcal{F} \geq \mathcal{V}_p^n \mathcal{G}$. $\mathcal{G} \xrightarrow{q} x$ implies there is $i \in I$, $x_i \in X_i$, and $\mathcal{H}_i \xrightarrow{q_i} x_i$ such that $\mathcal{G} \geq f_i(\mathcal{H}_i)$. Note that $\mathcal{F} \geq \mathcal{V}_p^n \mathcal{G} \geq \mathcal{V}_p^n f_i(\mathcal{H}_i) = f_i(\mathcal{V}_{p_i}^n \mathcal{H}_i)$. Since $\mathcal{H}_i \xrightarrow{q_i} x_i$, $\mathcal{V}_{p_i}^n (\mathcal{H}_i) \xrightarrow{\tau_{p_i} q_i} x_i$, and, thus, $\mathcal{F} \xrightarrow{s} x$.

The proof of (ii) is the similar. □

To avoid needless repetition, we state the next three corollaries to Theorem 4.6 only for the lower $p$-topological modifications. Analogous results obviously hold for the lower $p$-regular modifications as well.

**COROLLARY 4.7.** *Let $f:(X',q') \longrightarrow (X,q)$ be a convergence quotient map and $f:(X',p') \longrightarrow (X,p)$ an interior-preserving map. Then $f:(X',\tau_{p'} q') \longrightarrow (X,\tau_p q)$ is a convergence quotient map.*

**COROLLARY 4.8.** *Let $(X,q) = \sum_{i\in I}(X_i,q_i)$ be a disjoint sum of convergence spaces, and let $p \in \mathbb{C}(X)$ be such that each $g_i: (X_i,p_i) \longrightarrow (X,p)$ is an interior-preserving map, where $g_i: X_i \longrightarrow X$ is the canonical injection. Then $(X,\tau_p q) = \sum_{i\in I}(X_i,\tau_{p_i} q_i)$.*

**COROLLARY 4.9.** *Let $\Lambda = \{q_i: i \in I\} \subseteq \mathbb{C}(X)$ and let $p \in \mathbb{C}(X)$. If $q = \inf \Lambda$, then $\tau_p q = \inf \{\tau_p q_i: i \in I\}$.*

**THEOREM 4.10.** *Let $X$ be a set and let $\{(X_i,q_i): i \in I\}$ and $\{(X_i,p_i): i \in I\}$ be collec-*

tions of convergence spaces and, for all $i \in I$, let $f_i: X \to X_i$. Let $q$ be the initial structure on $X$ induced by $\{f_i: i \in I\}$ and $\{(X_i, q_i): i \in I\}$, and assume that $p$ is a structure such that $f_i: (X, p) \to (X_i, p_i)$ is continuous, for all $i \in I$.

(i) If $\tau^{p_i} q_i$ exists for all $i \in I$ and each $f_i: (X, p) \to (X_i, p_i)$ is an interior map, then $\tau^p q$ exists and is the initial structure on $X$ induced by $\{f_i: i \in I\}$ and $\{(X_i, \tau^{p_i} q_i): i \in I\}$.

(ii) If $r^{p_i} q_i$ exists for all $i \in I$ and each $f_i: (X, p) \to (X_i, p_i)$ is a closure map, then $r^p q$ exists and is the initial structure on $X$ induced by $\{f_i: i \in I\}$ and $\{(X_i, r^{p_i} q_i): i \in I\}$.

**PROOF.** (i) To show $\tau^p q$ exists, it suffices, by Theorem 4.3, to show that if $\mathcal{F} \geq \mathcal{V}_p^n(x)$ for some $x \in X$ and $n \in \mathbb{N}$, then $\mathcal{F} \xrightarrow{q} x$. Since $f: (X, p) \to (X_i, p_i)$ is continuous for all $i \in I$, $f_i(\mathcal{V}_p^n(x)) \geq \mathcal{V}_{p_i}^n(f_i(x))$, and since each $\tau^{p_i} q_i$ exists by assumption, $\mathcal{V}_{p_i}^n(f_i(x)) \xrightarrow{q_i} f_i(x)$ for each $i \in I$. Since $q$ is the initial structure, $\mathcal{V}_p^n(x) \xrightarrow{q} x$, and, hence, $\mathcal{F} \xrightarrow{q} x$. The remainder of the proof of (i) is straight-forward and is omitted.

The proof of (ii) exactly parallels that of (i). □

The corollaries of Theorem 4.10, like those of Theorem 4.6, are stated only for the upper *p*-topological modifications. The corresponding results involving upper *p*-regular modifications can be supplied by the reader.

**COROLLARY 4.11.** *Let $(X, q)$ be a subspace of $(X', q')$ and let $p' \in \mathbb{C}(X')$. Also let $(X, p)$ be a subspace of $(X', p')$ and assume that $\tau^{p'} q'$ exists. If $X$ is $p'$-open, then $(X, \tau^p q)$ is a subspace of $(X', \tau^{p'} q')$.*

**PROOF.** Since $X$ is $p'$-open in $(X', p')$, the identity map from $(X, p)$ into $(X', p')$ is a continuous interior map, and so, the conclusion follows from Theorem 4.10(ii). □

**COROLLARY 4.12.** *Let $(X, q) = \Pi_{i \in I}(X_i, q_i)$, let $p_i \in \mathbb{C}(X_i)$ be such that $\tau^{p_i} q_i$ exists for each $i \in I$, and let $p \in \mathbb{C}(X)$ be such that the $i$th projection map $\pi_i: (X, p) \to (X_i, p_i)$ is continuous interior map for all $i \in I$. Then $\tau^p q$ exists, and $(X, \tau^p q) = \Pi_{i \in I}(X_i, \tau^{p_i} q_i)$.*

**COROLLARY 4.13.** *Let $X$ be a set, $\Lambda = \{q_i: i \in I\} \subseteq \mathbb{C}(X)$, and let $p \in \mathbb{C}(X)$ be such that $\tau^p q_i$ exists for all $i \in I$. If $q = \sup \Lambda$, then $\tau^p q$ exists and $\tau^p q = \sup\{\tau^p q_i: i \in I\}$.*

**5. The topological series of a convergence space.** If $(X, q)$ is a convergence space, it is well known that there is a finest topology $\tau q$ coarser than $q$ and a finest regular convergence structure $rq$ coarser than $q$. These are the *topological* and *regular modifications* of $q$. However, neither $\tau q$-convergence nor $rq$-convergence can be described directly in terms of $q$-convergence. Consequently, descending ordinal series have been devised to "bridge the gap" between $q$ and these two lower modifications.

The *regularity series* $(r_\alpha q)$, introduced in [4] and studied also in [5], can be easily characterized by means of the lower *p*-regular modification for an arbitrary ordinal number $\alpha$ as follows: $r_\alpha q = r_{p_\alpha} q$, where $p_0 = \delta$, $p_1 = q$, $p_\alpha = r_{\alpha-1} q$ if $\alpha - 1$ exists, and $p_\alpha = \inf\{r_\beta q: \beta < \alpha\}$ if $\alpha$ is a limit ordinal. The least ordinal $\alpha$ for which $r_\alpha q = r_{\alpha+1} q$ is called the *length of the regularity series* and is denoted by $\ell_R q$. It is easy to verify that $r_\alpha q = rq$ if and only if $\alpha \geq \ell_R q$.

The *decomposition series* $(\pi_\alpha q)$, introduced in [3], is a descending ordinal sequence of pretopologies terminating in $\tau q$. Just as the regularity series gives an ordinal measure of how "non-regular" a given convergence space is, so likewise does the decomposition series measure how "non-topological" the given space is. However, the construction of the regularity and decomposition series are fundamentally so different that interactions or comparisons between them are difficult to find or interpret.

The existence of the lower $p$-topological modification and its dual relationship to the lower $p$-regular modification provide means for constructing a new descending ordinal sequence called the *topological series* $(\tau_\alpha q)$ of $(X,q)$ which, like the decomposition series, stretches between $q$ and $\tau q$. Following the preceding description of the regularity series, we define: $\tau_\alpha q = \tau_{p_\alpha} q$, where $p_0 = \delta, p_1 = q, p_\alpha = \tau_{\alpha-1} q$ if $\alpha - 1$ exists, and $p_\alpha = \inf \{\tau_{\beta q}: \beta < \alpha\}$ if $\alpha$ is a limit ordinal. The resulting *topological series* is the exact dual of the regularity series. It can be shown that the length of the topological series cannot exceed that of the decomposition series. Additional results pertaining to these and other related ordinal series will be published later.

## References

[1] C. H. Cook and H. R. Fischer, *Regular convergence spaces*, Math. Ann. **174** (1967), 1-7. MR 37#5837. Zbl 152.39603.

[2] H. J. Biesterfeldt Jr., *Regular convergence spaces*, Indag. Math. **28** (1966), 605-607. MR 34#5048. Zbl 178.56402.

[3] D. C. Kent and G. D. Richardson, *The decomposition series of a convergence space*, Czechoslovak Math. J. **23** (1973), no. 98, 437-446. MR 48 1134. Zbl 268.54002.

[4] ______, *The regularity series of a convergence space*, Bull. Austral. Math. Soc. **13** (1975), no. 1, 21-44. MR 51 13964. Zbl 308.54001.

[5] ______, *The regularity series of a convergence space. II*, Bull. Austral. Math. Soc. **15** (1976), no. 2, 223-243. MR 55 6343. Zbl 328.54001.

[6] ______, *p-regular convergence spaces*, Math. Nachr. **149** (1990), 215-222. MR 93a 54018. Zbl 724.54004.

[7] ______, *Convergence spaces and diagonal conditions*, Topology Appl. **70** (1996), no. 2-3, 167-174. CMP 96 14. Zbl 862.54002.

[8] H. J. Kowalsky, *Limesräume und Komplettierung*, Math. Nachr. **12** (1954), 301-340. MR 17,390b. Zbl 056.41403.

WILDE AND KENT: DEPARTMENT OF MATHEMATICS, WASHINGTON STATE UNIVERSITY, PULLMAN, WA 99164-3113, USA